\documentclass[preprint]{article}

\usepackage{amsmath}
\label{} 
\usepackage{amssymb}
\usepackage{amsthm}
\title{Abelian $1-$factorizations of complete multipartite graphs} 
\author{M. ~\textsc{Bogaerts} }
\newtheorem{theorem}{Theorem}

\newtheorem{proposition}[theorem]{Proposition}
\newtheorem{lemma}[theorem]{Lemma}
\newtheorem{problem}{Problem}
\newcommand{\ZZ}{\mathbb{Z}}
\begin{document}
\maketitle
\begin{abstract}
An automorphism group $G$ of a $1-$factorization of the complete multipartite graph $K_{m\times n}$ consists in permutations of the vertices of the graph mapping factors to factors. In this paper, we give a complete answer to the existence or non-existence problem of a $1-$factorization of $K_{m\times n}$ admitting an abelian group acting sharply transitively on  the vertices of the graph.
\end{abstract}

\section{Introduction}

\indent In this paper, we consider $1-$factorizations of complete multipartite graphs $K_{m\times n}$ (i.e. with $m$ parts consisting of $n$ vertices). In order to avoid confusion with the complete graph $K_m=K_{m\times 1}$, it is from now assumed that each part of $K_{m\times n}$ has at least two vertices (i.e. $n\ge 2$). An {\bf $r-$factor} of a graph $\Gamma$ is a spanning subgraph with all vertices of degree $r$. For $r=1$, the concept of $1-$factor coincide with the definition of a perfect matching. An {\bf $r-$factorization} of a graph $\Gamma$ is a decomposition of $\Gamma$ as a union of disjoint $r-$factors.\\
A group of permutations of the vertices of $\Gamma$ mapping factors to factors is called an {\bf automorphism group of the $r-$factorization}. If the action of an automorphism group is sharply transitive on the vertices, the factorization is said to be {\bf sharply transitive}.\\
In \cite{G23}, G. Mazzuoccolo and G. Rinaldi  consider the following problem:
\begin{problem}\label{question}
Given a finite group $G$ of even order, which graphs $\Gamma$ have a $1-$factorization admitting $G$ as an automorphism group with a sharply transitive action on the vertex set?
\end{problem}
The problem above is the natural generalization of the largely studied case in which we set $\Gamma$ to be the complete graph. The benchmark on this topic is the following theorem of Hartman and Rosa \cite{G6} from 1985:
\begin{theorem}\label{ExistCyc}
A complete graph $K_n$ admits a $1-$factorization with a cyclic automorphism group acting sharply transitively on the vertices if and only if $n$ is even and $n\neq 2^t$, $t \ge 3$.
\end{theorem}
Later on, other specified classes of groups have been considered by various authors; among the others we would like to recall the generalizations of Theorem \ref{ExistCyc} to the entire classes of finite abelian groups (Buratti, \cite{G7}) and finitely generated abelian groups (Bonvicini and Mazzuoccolo, \cite{G27}).
Even if Problem \ref{question} remains still open if $\Gamma$ is the complete graph, it is interesting to consider other possible choices for the graph $\Gamma$: for instance, the main result of \cite{G23} states the existence and non-existence of a $1-$factorization of the complete multipartite graph admitting a cyclic automorphism group acting sharply transitively on the vertex-set (a cyclic $1$-factorization from now on).
\begin{theorem}\cite{G28, G23}\label{ExistMR}
A cyclic $1-$factorization of $K_{m\times n}$:
\begin{itemize}
\item[$\bullet$] does not exist if $m \equiv 3 \mod 4$, $n=2d$ where $d$ is odd
\item[$\bullet$] does not exist if $m= 2^v d$ $n=2 d'$ where $d$ and $d'$ are odd and $v \ge 2$
\item[$\bullet$] does not exist if $m=p^v$, $n=2$, where $p$ is prime such that $p\equiv 1 \mod 4$  and $v\ge 1$, 
\item[$\bullet$] exists if $m= 2d$, $n=2 d'$ where $d$ and $d'$ are odd.
\item[$\bullet$] exists if $m= 2^vd$, $n=2^u d'$ where $d$ and $d'$ are odd and $u>1$. 
\item[$\bullet$] exists if $m= 2^vd$, $n=d'$ where $d$ and $d'$ are odd and $v\ge 1$. 
\item[$\bullet$] exists if $n=2$, $m\equiv 1 \mod 4 $ and  $m$ is not a prime power,
\item[$\bullet$] exists if $m\equiv 1 \mod 4 $ and $n= 2d$ where $d>1$ is odd.
\end{itemize}
\end{theorem} 
In this paper, we consider the following extension of problem \ref{question}:
\begin{problem}\label{question2}
Given a graph $\Gamma$, does a $1-$factorization of $\Gamma$ admitting an abelian group $G$ as automorphism group with a sharply transitive action on the vertex set exists?
\end{problem}
Of course, if a cyclic factorization exists, then this question has a positive answer. On the other hand, the trivial condition $mn\equiv 0 \mod 2$ leaves three main cases to consider, that are covered by our main theorem.\\
\begin{theorem}\label{ExistMB}
An abelian $1-$factorization of $K_{m\times n}$:
\begin{itemize}
\item[$\bullet$] exists if $m= 2^v d$ $n=2 d'$ where $d$ and $d'$ are odd and $v \ge 2$,
\item[$\bullet$] does not exist if $m=p$, $n=2$, where $p$ is prime such that $p\equiv 1 \mod 4$, 
\item[$\bullet$] exists if $m=p^v$, $n=2$, where $p$ is prime such that $p\equiv 1 \mod 4$  and $v\ge 2$, 
\item[$\bullet$] does not exist if $m \equiv 3 \mod 4$, $n=2d$ where $d$ is odd
\end{itemize}
\end{theorem}
The material involved in the constructions is presented in section \ref{Preliminaries}. Sections \ref{SecN2} to \ref{SecN4} state the existence and non-existence results. \\

\section{Preliminaries}\label{Preliminaries}

This paper is dedicated to the construction of $1-$facto\-ri\-zation of $K_{m\times n}$ ($mn$ even) with an abelian automorphism group $G$ acting sharply transitively on the vertices. For sake of clarity, the composition is written additively, with neutral element $\bar{0}$. As in \cite{G23}, the graph $\Gamma=K_{m\times n}$ is considered as the Cayley graph $Cay(G,\Omega)$ with $\Omega =G - H$ where $G$ is a group of order $mn$ and $H$ is a subgroup of order $n$ of $G$. Vertices of $\Gamma$ are the elements of $G$, and edges are (unordered) pairs of elements $[g_1,g_2], g_i\in G$ such that the difference $g_1-g_2$ belongs to $\Omega$. With this definition, the group $G$ acts naturally sharply transitively on its elements (i.e. on the vertices of $\Gamma$).\\

In \cite{G23}, G. Mazzuoccolo and G. Rinaldi give the definition of a starter for $K_{m\times n}$, a generalization of the concept defined by M. Buratti in \cite{G7}. The existence of a starter is proven to be equivalent to the existence of a $1-$facto\-ri\-zation of $K_{m\times n}$ with an automorphism group $G$ acting sharply transitively on the vertices. We here give the definition of a starter for $K_{m\times n}$ with a abelian group.\\

The set $\Omega$ can be partitioned as $\Omega_1\cup \Omega_2\cup \Omega_2'$ where $\Omega_1$ contains all involutions of $\Omega$, and for any $g\in \Omega_2$, $-g\in \Omega_2'$. The edge set of $\Gamma$ can be described as the union of orbits of $G$ acting additively, i.e. $(\cup_{j\in \Omega_1} \{[k, j+k]:k\in G\} ))\cup (\cup_{g\in \Omega_2} \{[l, l+g]: l \in G\}))$.  For any $j\in \Omega_1$, $Orb_G([\bar{0}, j])=\{[k, j+k]:k\in G\}$ has $\frac{mn}{2}$ elements and forms a $1-$factor of $\Gamma$. Because there are only $\frac{mn}{2}$ elements in each orbit, edges of this kind are called {\bf short edges}.\\

If $g\in \Omega_2$, then $Orb_G([\bar{0}, g])$ has $2n$ elements and is a union of cycles of $\Gamma$. These cycles can be written as $(x, x+g, x+2g,\dots,  x+(k-1)g)$ where $k$ is the order of $g$ and $x$ is a representative of one right coset of $<g>$ in $G$. The edges in these orbits are called {\bf long edges}.\\

Now we define two mappings $\partial$, $\phi$. The first one maps an edge to the differences of its vertices, and the second one maps an edge to its vertices. Both give one element if the edge is short, and $2$ if the edge is long.

$$\partial([x,y])=\left\{ \begin{array}{l l} \{x-y,y-x\}&\mbox{ if }[x,y]\mbox{ is long}\\
\{x-y\}&\mbox{ if }[x,y]\mbox{ is short}\\
\end{array}\right.$$

$$\phi([x,y])=\left\{ \begin{array}{l l} \{x,y\}&\mbox{ if }[x,y]\mbox{ is long}\\
\{x\}&\mbox{ if }[x,y]\mbox{ is short}\\
\end{array}\right.$$

Then for a set $S$, we define $\partial(S)=\cup_{e\in S} \partial(e)$  and $\phi(S)=\cup_{e\in S}\phi(e)$.\\
Now we give the definition of a {\bf starter} for the pair $(G, \Omega)$: it is a set $\Sigma=\{S_1,\dots , S_k\}$ of subsets of the edges $E(\Gamma)$, with $k$ subgroups $H_1,\dots,H_k$ of $G$ such that:
\begin{itemize}
\item[$\bullet$] the union of the differences $\partial S_1 \cup \dots\cup \partial S_k$ is $\Omega$, but without repetition (every element appears exactly once).
\item[$\bullet$]  for each $S_i$, $\phi(S_i)$ contains exactly one representative of each right coset of $H_i$ in $G$.
\item[$\bullet$] for any $H_i$ and any short edge $[x,y]\in S_i$, $H_i$ contains the involution $y-x=x-y$.
\end{itemize}

Each $H_i-$orbit of a set $S_i \in\Sigma$ is a $1-$factor of $\Gamma$ and the action of $G$ on the $1-$factors covers all edges of $\Gamma$ exactly once. By construction, the factorization is preserved under the action of $G$.\\

The construction of a starter can be simplified by the existence of a subgroup $A$ of $G$ of index $2$, as mentioned in \cite{G23}:

\begin{proposition}\cite{G23}\label{Index2Exist}  
Let $G$ be a finite group possessing a subgroup $A$ of index 2 and let $\Sigma'=\{S_1,\dots,S_t\}$ be a set of subsets of $E(\Gamma)$ together with subgroups $H_1,\dots,H_t$ which satisfy the second and third conditions of the definition of a starter. If $\partial S_1 \cup \dots\cup \partial S_t \subset A\cap \Omega$ and it does not contain repeated elements, then $\Sigma'$ can be completed to a starter for the pair $(G, \Omega)$.
\end{proposition}

\section{Abelian $1$-factorization of $K_{2^v d \times 2d'}$, $d, d'$ odd, $v\ge 2$}\label{SecN2}

First, consider the case $d'=1$. In \cite{G9}, Bonisoli and Labbate have proven the following proposition:
\begin{proposition}\label{ExistFixFact}
For $m=2^v d$, $v\ge 1$, $d$ odd, there exists a $1-$factorization of $K_{2m}$ with $C_m\times \ZZ_2$ as vertex-regular automorphism group, with one fixed factor $F$.
\end{proposition}

The graph obtained by deletion of the edges of $F$ is $K_{m\times 2}$, proving that this graph has a $1-$factorization preserved by a sharply transitive abelian automorphism group. We generalize this construction for other values as follows.\\

\begin{proposition}\label{ExistN2}
For $m=2^v d$, $n= 2d'$, $v\ge 2$, $d$ and $d'$ odd, there exists a $1-$factorization of $K_{m\times n}$ with $G= C_{d}\times C_{2^v d'} \times \ZZ_2$ as vertex-regular automorphism group.
\end{proposition}

The proof of this proposition relies on the following construction lemma.

\begin{lemma}\label{Doublinglemma} If there exists a starter for a cyclic regular $1-$factorization of $K_{m\times n}$ with automorphism group $G$, then there exists a starter for $K_{m\times 2n}$ with $G'=G\times \ZZ_2$.
\end{lemma}

{\bf Proof of lemmma \ref{Doublinglemma}}\\

By hypothesis, $K_{m\times n}$ has a starter for the cyclic group $G=\ZZ_{mn}$ with subgroup $H$ of order $n$. Let $\Sigma=\{S_1,\dots ,S_k\}$ together with subgroups $H_1, \dots, H_k$ be this starter. Now let consider $G'= G\times \ZZ_2$ with subgroup $H'= H\times \ZZ_2$, in such a way that $\Omega'=G' \backslash H'$ give $K_{m\times 2n}= Cay(G', \Omega')$.\\

Let  $\Sigma'$ denote the starter constructed with sets $\{S^1_1, \dots , S^1_k, S^2_1, \dots , S^2_k\}$ with subgroups  $H^1_1, \dots, H^1_k$, $ H^2_1, \dots, H^2_k$. For $i=1,\dots , k$, $S_i^1$ (resp. $S_i^2$) is defined as the set of pairs $ [(a,0), (b,0)]$ (resp. $ [(a,0), (b,1)]$) for all pairs $ [a,b] $ in $S_i$, and $H_i^1=H_i^2$ is defined as $H_i\times \ZZ_2$.\\

Elements in $\Omega'= G' \backslash H'$ can be written as $(a,c)$ with $a\in \Omega$ and $c\in \ZZ_2$. By construction, the union  $\delta S_1^1 \cup \dots \cup \delta S_k^1$ covers exactly once all elements of type $(a, 0)$ with $a\in \Omega$, while the union $\delta S_1^2 \cup \dots \cup \delta S_k^2$ does the same with all elements of type $(a, 1)$ with $a\in \Omega$. \hfill $\Box$\\

Since a cyclic regular $1-$factorization of $K_{2^v d\times d'}$ always exists, as stated in theorem \ref{ExistMR},  proposition \ref{ExistN2} is easily deduced from lemma \ref{Doublinglemma}.

\section{Abelian $1$-factorization of $K_{p^v  \times 2}$, $p\equiv 1 \mod 4$ prime, $v\ge 1$}\label{SecN3}

The case $K_{p^v d \times 2}$, $p\equiv 1 \mod 4$ prime, $v\ge 1$ is treated as follows. First, let remark that if $v=1$, up to isomorphism, the only possible abelian group is the cyclic group of order $2p$, for which it is known that a starter does not exists. For all $v\ge 2$, a starter always exists, as stated in the following proposition.\\
\begin{proposition}\label{ExistN3}
For $m=p^v$, $p\equiv 1 \mod 4$ prime , $v\ge 2$, there exists a $1-$factorization of $K_{m\times 2}$ with $G= \ZZ_{p^{v-1}}\times \ZZ_{p } \times \ZZ_2$ as vertex-regular automorphism group.
\end{proposition}

{\bf Proof}\\

The existence of the $1-$factorization will follow from the constuction of the starter.
First, the subgroup $H$ of $G= \ZZ_{p^{v-1}}\times \ZZ_{p } \times \ZZ_2$ is defined by its generator $(0,0,1)$. Let $p=4t+1$ and $t'=\frac{ p^{v-1} -1}{4}$, and let define the pairs $(S_i,H_i)$. \\

Let $H_1=<(1,0,0)>$. The set $S_1$ is defined as the union $(\cup_{j=1}^4 S_1^j)\cup \{e_1, e_2\}$ where:\\
$S_1^1 =\{ [(0,i, 0), (0, p-i, 0)] : i= 1, \dots, t\}$, \\
$S_1^2 =\{ [(0,i-1, 1), (0, p-i, 1)] : i= 1, \dots, t\}$,\\
$S_1^3 =\{ [(0,t+i, 0), (0, p-t-i, 1)] : i= 1, \dots, t\}$,\\
$S_1^4 =\{ [(0,t+i+1, 1), (0, p-t-i, 0)] : i= 1, \dots, t-1\}$,\\
$e_1= [(0,0,0), (2,t,1)]$, $e_2=[(0,2t+1,0), (2,t+1, 1)]$\\
%
%$\delta S_1^a=\{ (0, p-2i, 0), (0, 2i,0):  i= 1, \dots, t\}$\\
%$\delta S_1^b=\{ (0, p-2i+1, 0), (0, 2i-1,0):  i= 1, \dots, t\}$\\
%$\delta S_1^c =\{ [(0,p-2t-2i, 1), (0, 2t+2i, 1)] : i= 1, \dots, t\}$\\
%$\delta S_1^d =\{ [(0,p-2t-2i-1, 1), (0, 2t+2i+1, 1)] : i= 1, \dots, t-1\}$\\

$H_m$ is defined for $m=1,\dots , 2t'$ by: $H_m=H_1$,\\
$S_m^1=\{[(0, i, 0),(2m-1, p-i,0)]: i=1,\dots t  \} $,\\
$S_m^2=\{[(0, i-1, 1),(2m-1, p-i,1)]: i=1,\dots t $,\\
$S_m^3= \{[(0, p-t-i, 0),(2m, t+i,0)]i=1,\dots t \}$,\\
$S_m^4= \{[(0, p-t-i, 1),(2m, t+i-1,1)]i=1,\dots t \}$,\\
$S_m^5 = \{[(0,0,0),(2m-1,2t,1)]\}$\\
For each $m$, $S_m$ is defined by the union $\cup_{j=1}^5 S_m^j $\\

%
%$\delta S_m^a=\{ (2m-1, p-2i, 0), (p^{v-1}+1-2m, 2i,0):  i= 1, \dots, t\}$\\
%$\delta S_m^b=\{ (2m-1, p-2i+1, 0), (p^{v-1}+1-2m, 2i-1,0):  i= 1, \dots, t\}$\\
%$\delta S_m^c =\{ [(2m,p-2t-2i, 1), (p^{v-1}-2m, 2t+2i, 1)] : i= 1, \dots, t\}$\\
%$\delta S_m^d =\{ [(2m,p-2t-2i-1, 0), (p^{v-1}-2m, 2t+2i+1, 0)] : i= 1, \dots, t\}$\\
%$\delta S_m^e =\{ [(2m-1,2t, 1), (p^{v-1}-2m+1, p-2t, 1)] \}$\\
%
Now let $H_{2t'+1}= <(0,1,0)>$ and define $S_{2t'+1}:=\cup_{j=1}^4 S_{2t'+1}^j$ with the following subsets.\\
$S_{2t'+1}^1=\{[(1-i, 0,0), (i,0,0)]: i =1,\dots, t'\}$,\\
$S_{2t'+1}^2=\{[(i, 0,1), (-i,0,1)]: i =1,\dots, t'\}$,\\
$S_{2t'+1}^3=\{[(t'+i, 0,), (-t'-i,2t+2,1)]: i =1,\dots, 2t'\}$,\\
$S_{2t'+1}^4=\{[(0,2,1), (-t',0,0)]\}$\\
%
%
%$\delta S_{2t'+1}^a=\{ (p^{v-1}+1-2i, 0, 0), ( 2i-1,0,0):  i= 1, \dots, t'\}$\\
%$\delta S_{2t'+1}^b=\{ (p^{v-1}-2i,0, 0), (2i, 0,0):  i= 1, \dots, t'\}$\\
%$\delta S_{2t'+1}^c =\{ [(p^{v-1}-2t-2i,2t+2, 1), ( 2t+2i,p-2t-2, 1)] : i= 1, \dots, 2t'\}$\\
%$\delta S_{2t'+1}^d =\{ [(p^{v-1}-t', 2, 1),(t',p-2, 1) ] \}$\\
%
Let $A$ be the subgroup $\ZZ_{p^{v-1}}\times \ZZ_p$ of index $2$ of $G$. Let  $\Sigma'=\{S_1,\dots, S_{2t'+1} \}$ be the family of subsets of edges together with subgroups $(H_1,\dots, H_{2t'+1})$. It is a simple check that all elements $(i,j,k)$ of $\Omega$ with $k=0$ are covered by $\delta(\cup_l S_l)$, and that the second and third conditions of the definition of a starter are satisfied. Proposition \ref{Index2Exist} can be applied, resulting in the existence of the appropriate starter.   
\hfill $\Box$

\section{Abelian $1$-factorization of $K_{m\times 2d}$,  $m \equiv 3 \mod 4$,  $d$ odd}\label{SecN4}

\begin{proposition}\label{NExistN4}
For $m\equiv 3 \mod 4 $ , $d$ odd,  an abelian $1-$factorization of $K_{m\times 2d}$ with a sharply transitive action on the vertices does not exist.
\end{proposition}

{\bf Proof}\\

Here we use the well-known fact that any abelian group $G$ can be written as a direct product of cyclic group whose orders are equal to a power of a prime: $G=C_{p_1^{v_1}}\times \dots \times C_{p_k^{v_k}}$, $p_i$ prime for all $i=1,\dots, k$. Supose that there exists a starter for an abelian regular $1-$factorization of $K_{m\times 2d}$. The group $G$ must have order $2md$, and is isomorphic to the direct product $G'\times \ZZ_2$, where $G'$ the abelian subgroup of index $2$. The subgroup $H$ has order $2d$ and can be decomposed as $H=H'\times \ZZ_2$ where $H'$ is an abelian group order $d$. The set $\Omega= G\backslash H$ has $2md-2d=2d (m-1)$ elements, none of them is an involution. The number of elements in $\Omega$ of type $(a, 0)$ is equal  to the number of elements that can be written as $(a, 1)$. Since $\Omega$ contains no involutions, all edges are long. The subgroups $H_k$ choosen to construct the starter are of odd order, so the set $\phi (S_k)$ of representative of the cosets counts as many elements of type $(a, 0)$ than elements of type $(b, 1)$. A simple check reveals that the number of edges $e$ in $S_k$ with $\delta e = \{(a,0), (-a,0)\}$ is even. Each of these edges contributes to cover two elements of $\Omega$ with $\delta e$. For each $k$, $\delta S_k$ covers $4r$ elements of type $(a, 0)$ of $\Omega$. On the other hand, $\Omega$ contains $d(m-1)\equiv 2 \mod 4$ elements of this type, leading to a contradiction.
\hfill $\Box$\\

\bibliography{Notes}
\bibliographystyle{amsplain}

\end{document}